\documentclass[11pt,final,a4paper]{article}
\usepackage{amsmath,amssymb,setspace,amsthm}

\newcommand{\bm}[1]{\boldsymbol{#1}}

\newtheorem{theorem}{Theorem}[section]
\newtheorem{lemma}[theorem]{Lemma}      
\newtheorem{proposition}[theorem]{Proposition} 

\title{A polynomial generalization of the \\
power-compositions determinant\thanks{%
Work partially supported by the Ministerio de Ciencia y
Tecnolog\'ia under projects BFM2003-00368 and
 MTM2004-01728 and Ministerio de Ciencia y
Tecnolog\'ia and by the
   Generalitat de Catalunya under project 2005 SGR 00692} }
\author{
Josep M. Brunat and Antonio Montes \\Departament de Matem\`atica
Aplicada II, \\ Universitat Polit\`ecnica de Catalunya, Spain.\\
 \{josep.m.brunat, antonio.montes\}@upc.edu \\
http://www-ma2.upc.edu/$\sim$montes}

\date{July 12, 2005.}

\begin{document}

%
\setcounter{page}{1}

\maketitle

\begin{abstract}
  Let $C(n,p)$ be the set of $p$-compositions of an integer $n$, i.e., the set
  of $p$-tuples $\bm{\alpha}=(\alpha_1,\ldots,\alpha_p)$ of
  nonnegative integers such that $\alpha_1+\cdots+\alpha_p=n$, and
  $\mathbf{x}=(x_1,\ldots,x_p)$ a vector of indeterminates.  For
  $\bm{\alpha}$ and ${\bm{\beta}}$ two $p$-compositions of
  $n$, define $(\mathbf{x}+\bm{\alpha})^{\bm{\beta}}
  =(x_1+\alpha_1)^{\beta_1}\cdots (x_p+\alpha_p)^{\beta_p}$.  In this paper we
  prove an explicit formula for the determinant
  $\det_{\bm{\alpha},{\bm{\beta}}\in
    C(n,p)}((\mathbf{x}+\bm{\alpha})^{\bm{\beta}})$. In
  the case $x_1=\cdots=x_p$ the formula gives a proof of a conjecture
  by C.~Krattenthaler.
\bigskip


\end{abstract}

\section{Introduction}

Let us start with some notation. If
$\mathbf{u}=(u_1,\ldots,u_\ell)$ and $\mathbf{v}=(v_1,\ldots, v_{\ell})$ are
two vectors of the same length, we define
$\mathbf{u}^\mathbf{v}=u_1^{v_1}\cdots u_\ell^{v_\ell}$
(where, to be consistent $0^0=1$). In our case, the entries $u_i$ and $v_i$ of
$\mathbf{u}$ and $\mathbf{v}$ will be nonnegative integers or polynomials. We
use $\mathbf{x}=(x_1,\ldots,x_p)$ to denote a vector of indeterminates and
$\mathbf{1}=(1,\ldots,1)$.  The lengths of $\mathbf{x}$ and $\mathbf{1}$
will be clear from the context.  If $\mathbf{u}=(u_1,\ldots,u_\ell)$, then
$s(\mathbf{u})$ denotes the sum of the entries of $\mathbf{u}$, i.e.
$s(\mathbf{u})=u_1+\cdots+u_\ell$, and $\bar{\mathbf{u}}$ denotes the vector
obtained from $\mathbf{u}$ by deleting the last coordinate,
$\bar{\mathbf{u}}=(u_1,\ldots,u_{\ell-1})$.

Let $C(n,p)$ be the set of $p$-compositions of an integer $n$, i.e., the set
of $p$-tuples $\bm{\alpha}=(\alpha_1,\ldots,\alpha_p)$ of nonnegative
integers such that $\alpha_1+\cdots+\alpha_p=n$. If
$\bm{\alpha}=(\alpha_1,\ldots,\alpha_p)$ and
$\bm{\beta}=(\beta_1,\ldots,\beta_p)$ are two $p$-compositions of $n$,
using the above notation, we have
$\bm{\alpha}^{\bm{\beta}}=\alpha_1^{\beta_1}\cdots
\alpha_p^{\beta_p}$.  In~\cite{BrMo} the following explicit formula for the
determinant
$\Delta(n,p)=\det_{\bm{\alpha},\bm{\beta}\in C(n,p)}
\left({\bm{\alpha}}^{\bm{\beta}}\right)$  was proved:
\begin{equation}
\label{BM}
\Delta(n,p)=\prod_{k=1}^{\min\{n,p\}}\left( n^{n-1 \choose k}
    \prod_{i=1}^{n-k+1} i^{(n-i+1){n-i-1 \choose k-2}}\right)^{p\choose k}.
\end{equation}
In a complement~\cite{K2} to his impressive \emph{Advanced Determinant
Calculus}~\cite{K1}, C. Krattenthaler mentions this determinant, and after
giving the alternative formula
\begin{equation}
\label{K}
\Delta(n,p)=n^{n+p-1\choose p}\prod_{i=1}^{n} i^{(n-i+1){n+p-i-1\choose p-2}}
\end{equation}
he states as a conjecture a generalization to univariate polynomials. Namely,
let $x$ be an indeterminate and
$$
\Delta(n,p,x)=\det_{\bm{\alpha},\bm{\beta}\in C(n,p)}
\left((x\cdot \mathbf{1}+\bm{\alpha})^{\bm{\beta}}\right).
$$
Note that $(x\cdot \mathbf{1}+\bm{\alpha})^{\bm{\beta}}
=(x+\alpha_1)^{\beta_1}\cdots (x+\alpha_p)^{\beta_p}$.
\medskip

\noindent\textbf{Conjecture } [C. Krattenthaler]:
\begin{equation}
\label{Kxp}
\Delta(n,p,x)=
(px+n)^{{n+p-1 \choose p}}\prod_{i=1}^n i^{(n-i+1){n+p-i-1\choose p-2}}.
\end{equation}
As $(n-i+1){n+p-i-1\choose p-2}=(p-1){n+p-i-1\choose p-1}$,
formula~(\ref{K}) can be written in the form
$$
\Delta(n,p)=n^{n+p-1\choose p}\prod_{i=1}^{n} i^{(p-1){n+p-i-1\choose p-1}}
$$
and Krattenthaler's Conjecture~(\ref{Kxp}) in the form
\begin{equation}
\label{Kx}
\Delta(n,p,x)=
(px+n)^{{n+p-1 \choose p}}\prod_{i=1}^n i^{(p-1){n+p-i-1\choose p-1}}.
\end{equation}
The main goal of this paper is to prove a generalization of
formula~(\ref{Kx}) for $p$ indeterminates.
For this, let $\mathbf{x}=(x_1,\ldots,x_p)$ be a vector of
indeterminates, and let
$$
\Delta(n,p,\mathbf{x})=\det_{\bm{\alpha},\bm{\beta}\in
  C(n,p)}\left((\mathbf{x}+\bm{\alpha})^{\bm{\beta}}\right).
$$
(Recall that
$(\mathbf{x}+\bm{\alpha})^{\bm{\beta}}=(x_1+\alpha_1)^{\beta_1}\cdots
(x_p+\beta_p)^{\beta_p}$).  Then, we prove the following formula
(Theorem~\ref{TBMx}):
\begin{equation}
\label{BMx}
\Delta(n,p,\mathbf{x})
=(s(\mathbf{x})+n)^{n+p-1\choose p}
\prod_{i=1}^{n} i^{(p-1){n+p-i-1\choose p-1}}.
\end{equation}
As $s(\mathbf{x})=x_1+\cdots+x_p$, if $x_1=\cdots=x_p=x$, then
$s(\mathbf{x})=px$ and the conjectured identity~(\ref{Kx}) follows.

We also prove a variant of this result for proper compositions.  A
\emph{proper $p$-composition} of an integer $n$ is a $p$-composition
$\bm{\alpha}=(\alpha_1,\ldots,\alpha_p)$ of $n$ such that
$\alpha_i\ge 1$ for all $i=1,\ldots,n$. Denote by $C^*(n,p)$ the set of proper
$p$-compositions of $n$ and define
$$
\Delta^*(n,p,\mathbf{x})
=\det_{\bm{\alpha},\bm{\beta}\in C^*(n,p)}
\left((\mathbf{x}+\bm{\alpha})^{\bm{\beta}}\right).
$$
The determinant $\Delta^*(n,p,\mathbf{x})$ has the following
factorization (Theorem~\ref{TBMPx}):
\begin{equation}
\label{proper}
\Delta^*(n,p,\mathbf{x})=(s(\mathbf{x})+n)^{n-1\choose p}
\left(\prod_{j=1}^p
\prod_{i=1}^{n-p+1}(x_j+i)^{n-i-1\choose p-2}
\right)\prod_{i=1}^{n-p+1}i^{(p-1){n-i-1\choose p-1}}.
\end{equation}

The paper is organized as follows. In the next section we collect some
combinatorial identities for further reference.  In Section~\ref{equivalence}
we prove the equivalence between the formula~(\ref{K}) given by Krattenthaler
and~(\ref{BM}).  In Section~\ref{recurrence} we prove two lemmas. The first
one is a generalization of the determinant $\Delta(n,2,\mathbf{x})$. The
second lemma uses the first and corresponds to a property of a sequence of
rational functions which appear in the triangulation process of the
determinant $\Delta(n,p,\mathbf{x})$. Section~\ref{main} contains the proof
of the main result, Theorem~\ref{TBMx}.  Finally, Section~\ref{sproper} is
devoted to proving~(\ref{proper}).

\section{Auxiliary summation formulas}
\label{auxiliary}

\begin{lemma}
\label{ci}
Let $a,b,c,d,m$ and $n$ be nonnegative integers. Then, the
  following equalities hold.
\begin{enumerate}
\item[\rm (i)]
$\sum_{k\in\mathbb{Z}} {a\choose c+k}{b\choose d-k}
={a+b\choose c+d}$;
\item [\rm (ii)]
$\sum_{k\le n}{a+k\choose a}=\sum_{k\le n}{a+k\choose k}
={n+a+1\choose a+1}$;
\item[\rm (iii)]
$\sum_{r=1}^n r{n+a-r\choose a}={n+a+1\choose a+2}$;
\end{enumerate}
\end{lemma}
\begin{proof}%
(i) is the well known Vandermonde's convolution, see~\cite[p.~169]{GrKnPa}.
  The formulas in (ii) are versions of the parallel
  summation~\cite[p.~159]{GrKnPa}. Part (iii) follows from
{\footnotesize
\begin{eqnarray*}
\sum_{r=1}^n r{n+a-r\choose a}
&=&
\sum_{r=1}^n r{n+a-r\choose n-r}
=
\sum_{k=0}^{n-1} \sum_{i=0}^{k}{a+i\choose a} \\
&=& \sum_{k=0}^{n-1} {a+k+1\choose a+1}
= {a+n+1\choose a+2}.
\end{eqnarray*}
}
\end{proof}

\section{Equivalence between the two formulas for $\mathbf{x}=\mathbf{0}$}
\label{equivalence}
Here we prove the equivalence beetween the formulas~(\ref{BM}) and~(\ref{K})
for $\Delta(n,p)$. Obviously, the result of substituting $x=0$ in
formula~(\ref{Kxp}) of the Conjecture gives formula~(\ref{K}) for $\Delta(n,p)$.

\begin{proposition}
Formulas~\emph{(\ref{BM})} and~\emph{(\ref{K})} are equivalent.
\end{proposition}
\begin{proof}
  We derive formula~(\ref{K}) from~(\ref{BM}), which was
  already proved in~\cite{BrMo}. First, note that if $p<k\le n$, the binomial
  coefficient ${p\choose k}$ is zero. Thus, we can replace $\min\{p,n\}$ by $n$ in
  formula~(\ref{BM}). Analogously, if $n-k+1<i\le n$, the binomial
  coefficient ${n-i-1\choose k-2}$ is zero, and we can replace the upper value
  $n-k+1$ by $n$ in the inner product. Second, the case
  $a=n-1$, $b=d=p$ and $c=0$ of Lemma~\ref{ci} (i) yields
{\footnotesize
$$
\sum_{k=1}^n {n-1\choose k}{p\choose k} =
-1+\sum_{k=0}^n {n-1\choose k}{p\choose p-k}={n+p-1\choose p}-1,
$$
}%
and, if $i\ge 1$, by taking $a=n-i-1$, $b=d=p$ and $c=-2$
in Lemma~\ref{ci} (i), we obtain
{\footnotesize
$$
  \sum_{k=1}^{n-1}{n-i-1\choose k-2}{p\choose k}=
  \sum_{k}^{n-1}{n-i-1\choose k-2}{p\choose p-k}
  ={n+p-i-1\choose p-2}.
  $$
}
 Therefore,
\begin{eqnarray*}
\Delta(n,p)
&=& \prod_{k=1}^{\min\{n,p\}}\left( n^{n-1 \choose k}
    \prod_{i=1}^{n-k+1} i^{(n-i+1){n-i-1 \choose k-2}}\right)^{p\choose k}\\
&=&\prod_{k=1}^n\left( n^{n-1 \choose k}
    \prod_{i=1}^n i^{(n-i+1){n-i-1 \choose k-2}}\right)^{p
    \choose k}\\
&=& \left(\prod_{k=1}^n n^{{n-1\choose k}{p\choose k}}\right)
    \left(\prod_{k=1}^n
          \prod_{i=1}^n
          i^{(n-i+1){n-i-1\choose k-2}{p\choose k}}
    \right)\\
&=& n^{{n+p-1 \choose p}-1}
    \left(\prod_{i=1}^{n-1} i^{(n-i+1){n+p-i-1\choose p-2}}\right)
    n^{\sum_{k=1}^n {-1\choose k-2}{p\choose k}}\\
&=& n^{{n+p-1\choose p}+p-1}\prod_{i=1}^{n-1} i^{(n-i+1){n+p-i-1\choose p-2}}\\
&=& n^{n+p-1\choose p}\prod_{i=1}^{n} i^{(n-i+1){n+p-i-1\choose p-2}}.
\end{eqnarray*}
\end{proof}

\section{A recurrence}
\label{recurrence}
The next lemma evaluates the determinant
$$
D_r(n,y,z)=\det_{0\le i,j\le r}\left((y-i)^{n-j}(z+i)^j\right),
$$
by reducing it to a Vandermonde determinant. Note that
$D_n(n,x_1+n,x_2)=\Delta(n,2,\mathbf{x})$.
\begin{lemma}
\label{xy}
$$
D_r(n,y,z)=(y+z)^{r+1\choose 2}
\left(\prod_{i=0}^r (y-i)^{n-r}\right)
\left(\prod_{i=1}^r i^{r-i+1}\right).
$$
\end{lemma}
\begin{proof}
{\small
\begin{eqnarray*}
D_r(n,y,z)
&=&\left|\begin{array}{cccc}
(y-0)^n(z+0)^0 & (y-0)^{n-1}(z+0)^1 & \cdots & (y-0)^{n-r}(z+0)^r \\
(y-1)^n(z+1)^0 & (y-1)^{n-1}(z+1)^1 & \cdots & (y-1)^{n-r}(z+1)^r \\
\vdots         & \vdots             &        & \vdots             \\
(y-r)^n(z+r)^0 & (y-r)^{n-1}(z+r)^1 & \cdots & (y-r)^{n-r}(z+r)^r
\end{array}
\right| \\
&=&
\left( \prod_{i=0}^r (y-i)^n \right)
\left|\begin{array}{cccc}
1 & (z+0)/(y-0) & \cdots & (z+0)^r/(y-0)^r \\
1 & (z+1)/(y-1) & \cdots & (z+1)^r/(y-1)^r \\
\vdots & \vdots &        & \vdots          \\
1 & (z+r)/(y-r) &\cdots  & (z+r)^r/(y-r)^r
\end{array}
\right|  \\
&=&
\left(
\prod_{i=0}^r (y-i)^n\right)
\prod_{0\le i<j\le r}\left( \frac{z+j}{y-j}-\frac{z+i}{y-i} \right) \\
&=& \left(\prod_{i=0}^r (y-i)^n \right)
\prod_{0\le i<j\le r} \frac{(y+z)(j-i)}{(y-j)(y-i)} \\
&=&\left(
\prod_{i=0}^r (y-i)^n \right)
(y+z)^{r+1\choose 2}
\frac{\prod_{i=1}^r i^{r-i+1}}{\prod_{i=0}^r (y-i)^r} \\
&=&(y+z)^{r+1\choose 2}
\left(\prod_{i=0}^r (y-i)^{n-r}\right)
\left(\prod_{i=1}^r i^{r-i+1}\right).
\end{eqnarray*}
}
\end{proof}
\begin{lemma}
\label{rec}
Define
$f_r\colon \mathbb{N}_0\times\mathbb{N}_0 \rightarrow \mathbb{Q}(y,z)$
recursively by
\begin{eqnarray*}
f_0(i,j) &=& (z+i)^j; \\
f_{r+1}(i,j) &=& f_r(i,j) \quad \mbox{if}\quad j\le  r; \\
f_{r+1}(i,j) &=& f_r(i,j)
   -\left(\frac{y-i}{y-r}\right)^{j-r}\frac{f_r(i,r)f_r(r,j)}{f_r(r,r)}
\quad \mbox{if} \quad j> r.
\end{eqnarray*}
Then
\begin{enumerate}
\item[\rm (i)] $f_{r+1}(r,j)=0$ for $j\ge r+1$;
\item [\rm (ii)] $\displaystyle f_r(r,r)=(y+z)^r\frac{r!}{\prod_{i=0}^{r-1} (y-i)}$.
\end{enumerate}
\end{lemma}
\begin{proof}
  Part (i) is trivial using induction. To obtain $f_r=f_r(r,r)$, we take $n\ge
  r$ and calculate $D(n,y,z)=D_n(n,y,z)$ by Gauss triangulation method.

The entry $(i,j)$ of $D(n,y,z)$ is $(y-i)^{n-j}(z+i)^j=(y-i)^{n-j}f_0(i,j)$.
If $j\ge 1$, add to the column $j$ the column $0$ multiplied by
$$
-\frac{1}{(y-0)^{j-0}}\frac{f_0(0,j)}{f_0(0,0)}.
$$
Then, the entry $(i,j)$ with $j\ge 1$ is modified to
\begin{eqnarray*}
&&
(y-i)^{n-j}f_0(i,j)-(y-i)^{n-0}f_0(i,0)\frac{1}{(y-0)^{j-0}}\frac{f_0(0,j)}{f_0(0,0)} \\
&=&(y-i)^{n-j}\left\{ f_0(i,j)-\left(\frac{y-i}{y-0}\right)^{j-0}
\frac{f_0(i,k)f_0(k,j)}{f_0(0,0)}\right\} \\
&=& (y-i)^{n-j}f_{1}(i,j).
\end{eqnarray*}
Therefore,  $D(n,y,z)=\det_{0\le i,j\le r}\left((y-i)^{n-j}f_1(i,j)\right)$
and $f_1(0,j)=0$ for $j\ge 1$.

Now, assume that $D(n,y,z)=\det_{0\le i,j\le
  n}\left((y-i)^{n-j}f_k(i,j)\right)$ for $k\ge 1$ with $f_k(i,j)=0$ for
$k,j>i$. Add to the column $j\ge k+1$ the column $k$ multiplied by
$$
-\frac{1}{(y-k)^{j-k}}\frac{f_k(k,j)}{f_k(k,k)}.
$$
The entry $(i,j)$ is modified to
\begin{eqnarray*}
&&
(y-i)^{n-j}f_k(i,j)-(y-i)^{n-k}f_k(i,k)
\cdot
\frac{1}{(y-k)^{j-k}}\cdot\frac{f_k(k,j)}{f_k(k,k)}\\
&=&
(y-i)^{n-j}\left\{ f_k(i,j)-\left(\frac{y-i}{y-k}\right)^{j-k}
\frac{f_k(i,k)f_k(k,j)}{f_k(k,k)}\right\}\\
&=&
(y-i)^{n-j}f_{k+1}(i,j).
\end{eqnarray*}
Clearly $f_{k+1}(k,j)=0$ for $j>k$. After $n$ iterations, we get the
determinant of a triangular matrix. Hence
$$
D(n,y,z)=\det_{0\le k\le n}\left((y-k)^{n-k}f_k(k,k)\right)=\prod_{r=0}^n
(y-k)^{n-k}f_k.
$$
The principal  minor of order $r+1$ is
$D_r(n,y,z)=\prod_{k=0}^r (y-k)^{n-k}f_k$. Therefore,
\begin{equation}
\label{quotient}
\frac{D_r(n,y,z)}{D_{r-1}(n,y,z)}=(y-r)^{n-r}f_r.
\end{equation}
On the other hand, by Lemma~\ref{xy} we obtain
\begin{eqnarray*}
\frac{D_r(n,y,z)}{D_{r-1}(n,y,z)}
&=&
\frac{ (y+z)^{r+1\choose 2}
       \left(\prod_{i=0}^r (y-i)^{n-r}\right)
       \left(\prod_{i=1}^r i^{r-i+1}\right)  }
     { (y+z)^{r\choose 2}
       \left(\prod_{i=0}^{r-1} (y-i)^{n-r-1} \right)
       \left(\prod_{i=1}^{r-1} i^{r-i} \right)} \\
&=&
(y+z)^r\cdot r!\cdot \frac{(y-r)^{n-r}}{\prod_{i=0}^{r-1} (y-i)}.
\end{eqnarray*}
Comparing with~(\ref{quotient}), we have arrived at
$$
f_r=(y+z)^r\frac{r!}{\prod_{i=0}^{r-1} (y-i)}.
$$
\end{proof}

\section{Proof of the main theorem}
\label{main}
We sort $C(n,p)$ in lexicographic order. For instance, for $n=5$, and
$p=3$, we obtain
$$
\begin{array}{rl}
C(5,3)=& \{ \,         (5,0,0), (4,1,0), (3,2,0), (2,3,0), (1,4,0), (0,5,0), \\
       &\phantom{\{}\, (4,0,1), (3,1,1), (2,2,1), (1,3,1), (0,4,1), \\
       &\phantom{\{}\, (3,0,2), (2,1,2), (1,2,2), (0,3,2), \\
       &\phantom{\{}\, (2,0,3), (1,1,3), (0,2,3), \\
       &\phantom{\{}\, (1,0,4), (0,1,4), \\
       &\phantom{\{}\, (0,0,5) \, \}.
\end{array}
$$
Let $M(n,p,\mathbf{x})$ be the matrix with rows and columns labeled by the
$p$-compositions of $n$ in lexicographic order and with the entry
$(\bm{\alpha},\bm{\beta})$ equal to
$(\mathbf{x}+\bm{\alpha})^{\bm{\beta}}$.
We have $\Delta(n,p,\mathbf{x})=\det M(n,p,\mathbf{x})$.

An entry $(\mathbf{x}+\bm{\alpha})^{\bm{\beta}}$ in
$M(n,p,\mathbf{x})$ can be written in the form
$(\bar{\mathbf{x}}+\bar{\bm{\alpha}})^{\bar{\bm{\beta}}}
(x_p+\alpha_p)^{\beta_p}$. For $0\le i,j\le n$, let $S_{ij}$ be the matrix
with entries $(\bar{\mathbf{x}}+\bar{\bm{\alpha}})^{\bar{\bm{\beta}}}$ where $\bm{\alpha}$ and
$\bm{\beta}$ satisfy $\alpha_p=i$ and $\beta_p=j$.  Thus, the submatrix of
$M(n,p,\mathbf{x})$ formed by the entries labeled $(\bm{\alpha},\bm{\beta})$ with
$\alpha_p=i$ and $\beta_p=j$ can be written $\left(S_{ij}(x_p+i)^j\right)$.
Note that
$$
S_{kk}=M(n-k,p-1,\bar{\mathbf{x}}).
$$
Define $f_0(i,j)=(x_p+i)^j$. Therefore, $M(n,p,\mathbf{x})$ admits the block
decomposition
$$
M(n,p,\mathbf{x})=(S_{ij}f_0(i,j))_{0\le i,j\le n}.
$$
The idea is to put $M(n,p,\mathbf{x})$ in block triangular form
in such a way that at each step only the last factor of each block is
modified.

\begin{theorem}
\label{TBMx}
$$
\Delta(n,p,\mathbf{x})=
(s(\mathbf{x})+n)^{{n+p-1 \choose p}}
\prod_{i=1}^n i^{(p-1){n+p-i-1\choose p-1}}.
$$
\end{theorem}
\begin{proof}
The proof is by induction on $p$. For $p=1$, $\Delta(n,p,x)$ is the
determinant of the $1\times 1$ matrix $((x+n)^n)$. Hence
$\Delta(n,p,x)=(x+n)^n$. This value coincides with the right hand side of the
formula for $p=1$.

Consider now the case $p=2$. Any 2-composition of $n$ is of the form $(n-i,i)$
for some $i$, $0\le i\le n$.
The determinant to be calculated is
$\Delta(n,2,\mathbf{x})
= \det_{0\le i,j\le n} \left((x_1+n-i)^{n-j}(x_2+i)^j\right)$.
By taking $r=n$, $y=x_1+n$ and $z=x_2$ in Lemma~\ref{xy}, we get
$$
\Delta(n,2,\mathbf{x})=D_n(n,x_1+n,x_2)=
(x_1+x_2+n)^{n+1\choose 2}\prod_{i=1}^n i^{n-i+1}.
$$
Therefore, the formula holds for $p=2$.

Now, let $p>2$ and assume that the formula holds for $p-1$. Begin with the
block decomposition of the matrix $M(n,p,\mathbf{x})=(S_{ij}f_0(i,j))_{0\le
  i,j\le n}$.

Assume $\Delta(n,p,\mathbf{x})=\det (S_{ij}f_r(i,j))$ where
$S_{ij}
=((\bar{\mathbf{x}}+\bar{\bm{\alpha}})^{\bar{\bm{\beta}}})$,
\ with $\alpha_p=i$, \ $\beta_p=j$, and $f_r(i,j)=0$ for $i<r$ and $j>i$.

Fix a column $\bm{\beta}$ with $\beta_p=j>r$.  For each
$\bm{\gamma}\in C(n,p)$ with $\gamma_p=r$ and $\gamma_k\ge \beta_k$
for $k\in[p-1]$, add to the column $\bm{\beta}$ the column
$\bm{\gamma}$ multiplied by
$$
-\frac{1}{(s(\bar{\mathbf{x}})+n-r)^{j-r}}
{j-r \choose \bar{\bm{\gamma}}-\bar{\bm{\beta}}}
\frac{f_r(r,j)}{f_r(r,r)}.
$$
The differences
$\bar{\bm{\delta}}=\bar{\bm{\gamma}}-\bar{\bm{\beta}}$
are exactly the $(p-1)$-compositions of $j-r$. Also note that by the multinomial
theorem,%
{\small
$$
\left(s(\bar{\mathbf{x}})+n-i\right)^{j-r}
=\left((x_1+\alpha_1)+\cdots+(x_{p-1}+\alpha_{p-1})\right)^{j-r}
=\sum_{\bar{\bm{\delta}}}{j-r \choose \bar{\bm{\delta}}}
\left(s(\bar{\mathbf{x}})+\bar{\bm{\alpha}}\right)^{\bar{\delta}}.
$$
}%
Then, a term of column $\bm{\beta}$ is modified to
{\footnotesize
\begin{eqnarray*}
&&
(\bar{\mathbf{x}}+\bar{\bm{\alpha}})^{\bar{\bm{\beta}}}f_r(i,j)
-\sum_{\bar{\bm{\gamma}}}\frac{1}{(s(\bar{\mathbf{x}})+n-r)^{j-r}}
{j-r \choose \bar{\bm{\gamma}}-\bar{\bm{\beta}}}
\frac{f_r(r,j)}{f_r(r,r)}(\bar{\mathbf{x}}
+\bar{\bm{\alpha}})^{\bar{\bm{\gamma}}}f_r(i,r) \\
&=&
(\bar{\mathbf{x}}+\bar{\bm{\alpha}})^{\bar{\bm{\beta}}}
\left\{
f_r(i,j)-\frac{1}{(s(\bar{\mathbf{x}})+n-r)^{j-r}}
\left(
\sum_{\bar{\bm{\delta}}}{j-r \choose \bar{\bm{\delta}}}
(\bar{\mathbf{x}}+\bar{\bm{\alpha}})^{\bar{\bm{\delta}}}
\right)
\frac{f_r(r,j)f_r(i,r)}{f_r(r,r)}
\right\}\\
&=&
(\bar{\mathbf{x}}+\bar{\bm{\alpha}})^{\bar{\bm{\beta}}}
\left\{
f_r(i,j)-\frac{(s(\bar{\mathbf{x}})+n-i)^{j-r}}{(s(\bar{\mathbf{x}})+n-r)^{j-r}}
\frac{f_r(r,j)f_r(i,r)}{f_r(r,r)}
\right\}.
\end{eqnarray*}
} 
Now, define
$f_{r+1}(i,j)=f_r(i,j)$ for $j\le r$ and
$$
f_{r+1}(i,j)
=f_r(i,j)-\frac{(s(\bar{\mathbf{x}})+n-i)^{j-r}}{(s(\bar{\mathbf{x}})+n-r)^{j-r}}
\frac{f_r(r,j)f_r(i,r)}{f_r(r,r)}
$$
for $j>r$. Note that $f_{r+1}(r,j)=0$ for $j>r$. After $n$ iterations, we
arrive at the block matrix $(S_{ij}f_n(i,j))_{0\le i,j\le n}$ where $f(i,j)=0$ for $j>i$. Thus, the
determinant $\Delta(n,p,\mathbf{x})$ is the product of the determinants of the
diagonal blocks:
$$
\Delta(n,p,\mathbf{x})=\prod_{r=0}^n \det (S_{rr}f_r(r,r)).
$$
Now, $S_{rr}=M(n-r,p-1,\bar{\mathbf{x}})$, a square matrix of
order ${n-r+p-2\choose p-2}$. Therefore
$$
\Delta(n,p,\mathbf{x})=\prod_{r=0}^n
\left(\Delta(n-r,p-1,\bar{\mathbf{x}}) f_r(r,r)^{{n-r+p-2\choose p-2}}\right).
$$
Now, observe that the rational funcions $f_r$ satisfy the hypothesis of
Lema~\ref{rec} with $y=s(\bar{\mathbf{x}})+n=x_1+\cdots+x_{p-1}+n$ and
$z=x_p$. Thus,
$$
f_r=f_r(r,r)=(s(\mathbf{x})+n)^r\cdot \frac{r!}{\prod_{i=0}^{r-1}
  (s(\bar{\mathbf{x}})+n-i)}.
$$
By the induction hypothesis,
\begin{eqnarray*}
\Delta(n,p,\mathbf{x})
&=&
\prod_{r=0}^n \left(
(s(\bar{\mathbf{x}})+n-r)^{n-r+p-2\choose p-1}
\prod_{i=1}^{n-r} i^{(p-2){n-r+p-i-2\choose p-2}}
\right) \\
&& \cdot
\prod_{r=0}^n \left(
(s(\mathbf{x})+n)^r\cdot r!\cdot \frac{1}
{\prod_{i=0}^{r-1}(s(\bar{\mathbf{x}})+n-i)}
\right)^{n-r+p-2\choose p-2}
\end{eqnarray*}
It remains to count how many factors of each type there are in the above
product.

The number of factors $(s(\mathbf{x})+n)$ is
$\sum_{r=1}^{n} r{n+p-r-2\choose p-2}$.
From Lemma~\ref{ci} (iii) for $a=p-2$ this coefficient is
${n+p-1\choose p}$.

The number of factors $s(\bar{\mathbf{x}})+n-i$, for $0\le i\le n-1$, is
(by using Lemma~\ref{ci} (ii) with $a=p-2$)
{\footnotesize
$$
{n-i+p-2 \choose p-1}-\sum_{r=i+1}^{n}{n-r+p-2\choose p-2}
=
{n-i+p-2\choose p-1}-{n-i+p-2\choose p-1}
= 0.
$$}
Finally, for $1\le i\le n$, the number of factors equal to $i$ is
{\footnotesize
\begin{eqnarray*}
\lefteqn{(p-2)\sum_{r=0}^{n-i}{n+p-i-r-2\choose p-2}
+\sum_{r=i}^n {n+p-r-2\choose p-2}=} \\
&&
(p-2){n+p-i-r-1\choose p-1}+{n+p-r-1\choose p-1}=
(p-1){n+p-r-1\choose p-1}.
\end{eqnarray*}
}%
\end{proof}

\section{Proper compositions}
\label{sproper}
A \emph{proper $p$-composition} of an integer $n$ is a $p$-composition
$\bm{\alpha}=(\alpha_1,\ldots,\alpha_p)$ of $n$ such that $\alpha_i\ge
1$ for all $i=1,\ldots,n$.  We denote by $C^*(n,p)$ the set of proper
$p$-compositions of $n$.  In~\cite{BrMo} the following formula was given:
$$
\Delta^*(n,p)=\det_{\bm{\alpha},\bm{\beta}\in C^*(n,p)}
(\bm{\alpha}^{\bm{\beta}}) =n^{n-1\choose
  p}\prod_{i=1}^{n-p+1} i^{(n-i+1){n-i-1\choose p-2}}.
$$
Here, we study the corresponding generalization
$$
\Delta^*(n,p,\mathbf{x})=\det_{\bm{\alpha},\bm{\beta}\in C^*(n,p)}
\left((\mathbf{x}+\bm{\alpha})^{\bm{\beta}}\right).
$$

\begin{theorem}
\label{TBMPx} If $p\le n$, then
$$
\Delta^*(n,p,\mathbf{x})=
 (s(\mathbf{x})+n)^{n-1\choose p}
\left(\prod_{i=1}^{n-p+1}\prod_{j=1}^p (x_j+i)^{n-i-1\choose p-2}
\right)
\prod_{i=1}^{n-p}i^{(p-1){n-i-1\choose p-1}}.
$$
\end{theorem}
\begin{proof}
The mapping $C^*(n,p)\rightarrow C(n-p,p)$ defined by
$\bm{\alpha}=(\alpha_1,\ldots,\alpha_p)
\mapsto \bm{\alpha}-\mathbf{1}
=(\alpha_1-1,\ldots, \alpha_p-1)$
is bijective. Thus, we have
\begin{eqnarray*}
\Delta^*(n,p,\mathbf{x})
&=&
\det_{\bm{\alpha},\bm{\beta}\in C^*(n,p)}
\left((\mathbf{x}+\bm{\alpha})^{\bm{\beta}}\right) \\
&=&
\det_{\bm{\alpha},\bm{\beta}\in C^*(n,p)}
\left(
(\mathbf{x}+\mathbf{1}+\bm{\alpha}-\mathbf{1})^{\bm{\beta}-\mathbf{1}+\mathbf{1}}
\right) \\
&=&
\det_{\bm{\alpha},\bm{\beta}\in C(n-p,p)}
\left(
(\mathbf{x}+\mathbf{1}+\bm{\alpha})^{\bm{\beta}}
(\mathbf{x}+\mathbf{1}+\bm{\alpha})^\mathbf{1}
\right)\\
&=&
\Delta(n-p,p,\mathbf{x}+\mathbf{1})
\prod_{\bm{\alpha}\in C(n-p,p)}
(\mathbf{x}+\mathbf{1}+\bm{\alpha})^\mathbf{1}.
\end{eqnarray*}
The number of times that an integer $i$, $0\le i\le n-p$ appears as the first
entry of $p$-compositions of $n-p$ is the number of solutions
$(\alpha_2,\ldots, \alpha_{n-p})$ of $i+\alpha_2+\cdots+\alpha_p=n-p$, which
is ${n-p-i+p-2\choose p-2}={n-i-2\choose p-2}$. The count is the same for
every coordinate. Then, in the product $\prod_{\bm{\alpha}\in
  C(n-p,p)}(\mathbf{x}+\mathbf{1}+\bm{\alpha})^\mathbf{1}$, the number
of factors equal to $x_j+1+i$ is ${n-i-2\choose p-2}$; equivalently, for $1\le
i\le n-p+1$, the number of factors equal to $x_j+i$ is ${n-i-1\choose p-2}$.
Therefore, {\small
\begin{eqnarray*}
\Delta^*(n,p,\mathbf{x})&=&\Delta(n-p,p,\mathbf{x}+\mathbf{1})
     \prod_{\bm{\alpha}\in C(n-p,p)}
       (\mathbf{x}+\mathbf{1}+\bm{\alpha})^\mathbf{1}\\
&=& (s(\mathbf{x})+n)^{n-1\choose p}
\left(\prod_{i=1}^{n-p+1}\prod_{j=1}^p (x_j+i)^{n-i-1\choose p-2}
\right)
\prod_{i=1}^{n-p}i^{(p-1){n-i-1\choose p-1}}.
\end{eqnarray*}
}
\end{proof}

\section*{Acknowledgements}
The authors would like to thank the referee very much for valuable
suggestions, corrections and comments, which result in a great improvement of
the original manuscript.

\end{document}